\documentclass[letterpaper, 10 pt, conference]{ieeeconf}  

\IEEEoverridecommandlockouts                              
\usepackage{cite}
\usepackage{amsmath,amssymb,amsfonts,bm}
\usepackage{bbm}
\usepackage{graphicx}
\usepackage{subfigure}
\usepackage{textcomp}
\usepackage{algorithm}
\usepackage{algpseudocode}

\usepackage{xcolor}
\allowdisplaybreaks[4]
\newtheorem{definition}{Definition}
\newtheorem{theorem}{Theorem}
\newtheorem{lemma}{Lemma}

\title{\LARGE \bf
Data-enabled Policy Optimization for the Linear Quadratic Regulator
}

\author{Feiran Zhao, Florian D\"{o}rfler, Keyou You
\thanks{Research of F. Zhao and K. You was supported by National Key R\&D Program of China (2022ZD0116700) and National Natural Science Foundation of China (62033006, 62325305).}
\thanks{F. Zhao and K. You  are with the Department of Automation and BNRist, Tsinghua University, Beijing 100084, China. (e-mail: zhaofr18@mails.tsinghua.edu.cn, youky@tsinghua.edu.cn.) F. D\"{o}rfler is with the Department of Information Technology and Electrical Engineering, ETH Zurich, 8092 Zurich, Switzerland. (e-mail: dorfler@control.ee.ethz.ch)}}%
\begin{document}

\maketitle
\thispagestyle{empty}
\pagestyle{empty}

\begin{abstract}
Policy optimization (PO), an essential approach of reinforcement learning for a broad range of system classes, requires significantly more system data than indirect (identification-followed-by-control) methods or behavioral-based direct methods even in the simplest linear quadratic regulator (LQR) problem. In this paper, we take an initial step towards bridging this gap by proposing the data-enabled policy optimization (DeePO) method, which requires only a finite number of sufficiently exciting data to iteratively solve the LQR problem via PO. Based on a data-driven closed-loop parameterization, we are able to directly compute the policy gradient from a batch of persistently exciting data. Next, we show that the nonconvex PO problem satisfies a projected gradient dominance property by relating it to an equivalent convex program, leading to the global convergence of DeePO. Moreover, we apply regularization methods to enhance certainty-equivalence and robustness of the resulting controller and show an implicit regularization property. Finally, we perform simulations to validate our results.
\end{abstract}


\section{Introduction}

As a cornerstone of modern control theory, the linear quadratic regulator (LQR) problem has been the benchmark for data-driven control methods that seek to design a controller from raw system data. The manifold approaches to data-driven control can be broadly categorized as \textit{indirect} (when identifying a dynamical model followed by model-based control design) versus \textit{direct} (when bypassing the identification step). The use of direct data-driven control is usually motivated when the dynamical model is difficult to establish, or is too complex for model-based control design. As an end-to-end approach, the direct methods are conceptually simple and easy to implement in practice.

A representative instance of direct data-driven control is policy optimization (PO), an essential approach for applications of reinforcement learning (RL)~\cite{mnih2015human-level,lillicrap2016continuous,recht2019tour}. As an iterative method, PO directly searches over the policy space to optimize a performance metric of interest. Based on zeroth-order optimization techniques, it uses multiple system trajectories to estimate the policy gradient. There has been a resurgent interest in studying theoretical properties of PO on the LQR problem such as convergence and sample complexity; see e.g., \cite{fazel2018global, mohammadi2022convergence, malik2019derivative, zhao23global} and the comprehensive survey~\cite{bin2022towards}. Even though global convergence has been shown for the nonconvex PO problem by a \textit{gradient dominance} property~\cite{fazel2018global}, there exists a considerable gap in the sample complexity between PO and indirect methods, which have proved themselves to be more sample-efficient~\cite{tu2019gap,simchowitz2020naive} for solving the LQR problem. This gap is due to the exploration or trial-and-error nature of RL, or more specifically, that the cost used for gradient estimate can only be evaluated \textit{after} a whole trajectory is observed. Thus, the existing PO methods require numerous system trajectories to find an optimal policy, even in the simplest LQR setting.

Recent years have witnessed an emerging line of direct methods inspired by the \textit{Fundamental Lemma}~\cite{willems2005note}, which states that the behavior of a linear time-invariant (LTI) system can be characterized by the range space of raw data matrices. This result implies a non-parametric representation of LTI systems, giving rise to a notable implicit design called data-enabled predictive control (DeePC) \cite{coulson2019data}, which has seen many successful implementations in different practical scenarios~\cite{markovsky2021behavioral}. The fundamental lemma has also been utilized to solve various explicit control design and analysis problems~\cite{de2019formulas,van2020data,berberich2020data}. In particular, it has been shown in \cite{de2019formulas} that using subspace relations, the closed-loop LTI system can be parameterized by input-state data, leading to a data-based convex reformulation of the LQR problem. Compared with PO, this approach is significantly more sample-efficient as it only requires a batch of persistently exciting (PE) data. Indeed, the PE condition is equivalent to identifiability for LTI systems and should be a minimal assumption for most control design problems~\cite{van2020data,kang2022minimum}, e.g., the LQR problem. There have been many recent works leveraging regularization methods to promote certainty-equivalence and robustness of the LQR~\cite{de2021low,dorfler2021certainty,dorfler22on}, and to bridge behavioral-based direct and indirect methods~\cite{dorfler21bridging}. All these methods use only a small batch of PE data compared to data-hungry zeroth-order PO methods~\cite{fazel2018global, mohammadi2022convergence, malik2019derivative}. This leads to a natural question: does there exist a data-efficient PO method for solving the LQR problem?

In this paper, we provide an affirmative answer to the above question. By leveraging the data-driven closed-loop parameterization~\cite{de2019formulas}, we propose an iterative method called \textbf{d}ata-\textbf{e}nabl\textbf{e}d \textbf{p}olicy \textbf{o}ptimization (DeePO) to solve the LQR problem. Instead of estimating the policy gradient from the cost of observed trajectories, we show that after a change of optimization variables, the gradient can be directly characterized from a batch of PE data. Even though the resulting optimization problem is nonconvex, it can be parameterized as a data-based convex program. By exploiting this relation and using a recent PO result~\cite{sun2021learning}, we further show that the LQR cost is \textit{projected gradient} dominated, while it is only \textit{gradient} dominated in \cite{fazel2018global, mohammadi2022convergence}. By establishing that the cost is also locally smooth, we show that the projected gradient method converges to the global optimum. We also investigate how regularization~\cite{de2021low,dorfler2021certainty,dorfler22on} affects the convergence of DeePO. In particular, we show that the certainty-equivalence regularizer leads to an \textit{implicit regularization} property, meaning that the DeePO algorithm without regularization behaves as if it is regularized. This property has been advocated as an important feature of gradient-based methods for solving many nonconvex problems~\cite{neyshabur2017geometry,arora2019implicit,zhang2021policy}. Finally, we perform a numerical case study to validate our theoretical results. We are hopeful that the discovered DeePO method with significantly relaxed data requirements offers a possible path towards direct adaptive LQR control.

The rest of this paper is organized as follows. In Section \ref{sec:prob}, we revisit the LQR problem and recapitulate the data-driven LQR formulation. In Section \ref{sec:deepo}, we propose the DeePO method to iteratively solve the LQR problem and show its global convergence. Section \ref{sec:reg} studies the effects of two regularizers on the convergence of DeePO. Section \ref{sec:simu} uses a numerical example to validate our main results. Conclusion and future work in Section \ref{sec:conc} complete this paper.

\textbf{Notation.} We use $I_n$ to denote the $n$-by-$n$ identity matrix. We use $\underline{\sigma}(\cdot)$ to denote the minimal singular value of a matrix. We use $\|\cdot\|$ to denote the $2$-norm of a vector or a matrix, and $\|\cdot\|_F$ the Frobenius norm. We use $\rho(\cdot)$ to denote the spectral radius of a square matrix. We use $\text{poly}(\cdot)$ to denote a polynomial function. We use $\dagger$ to denote the right inverse of a full row rank matrix. 


\section{Problem Formulation}\label{sec:prob}

In this section, we first revisit the model-based LQR problem. By recapitulating its direct data-driven formulation from \cite{de2019formulas}, we then propose our PO reformulation.

\subsection{The Model-based LQR problem}
Consider a discrete-time LTI system
\begin{equation}\label{equ:sys}
x(t+1) = Ax(t) + Bu(t),
\end{equation}
where $x(t)\in\mathbb{R}^{n}$ and $u(t)\in\mathbb{R}^{m}$ are the state and control input, respectively. We assume that $(A,B)$ are controllable.

The LQR problem is phrased as finding a state-feedback gain $K\in \mathbb{R}^{m\times n}$ to minimize the quadratic cost
\begin{equation}\label{prob:lqr} J(K):=\mathbb{E}_{x(0)\sim \mathcal{D}}\left[\sum_{t=0}^{\infty}(x(t)^{\top} Q x(t)+u(t)^{\top} R u(t))\right],
\end{equation}
where $Q\succ 0, R \succ 0$ are penalty matrices, and $\{x(t),u(t)\}$ is the trajectory following (\ref{equ:sys}) and $u(t)=Kx(t)$ starting from the initial state $x(0)$. The distribution $\mathcal{D}$ of $x(0)$ satisfies $\mathbb{E}[x(0)]=0$ and $\mathbb{E}[x(0)x(0)^{\top}]=I_n$. It is well-known that the unique optimal gain to (\ref{prob:lqr}) is
$$
	K^* = - (R+B^{\top}P^*B)^{-1}B^{\top}P^*A,
$$
where $P^*$ is the unique positive semi-definite solution to the algebraic Riccati equation~\cite{bertsekas1995dynamic}
$$
P^*=A^{\top} P^* A+Q-A^{\top} P^* B(R+B^{\top} P^* B)^{-1} B^{\top} P^* A.
$$

We aim to solve the LQR problem in a direct data-driven approach when $(A,B)$ are unknown, but we assume the access to a $T$-length dataset of states and control inputs.

\subsection{Direct data-driven formulation}
Define the offline data matrices
\begin{align*}
X_- &= \begin{bmatrix}
x(0)& x(1)& \dots& x(T-1)
\end{bmatrix}\in \mathbb{R}^{n\times T},\\
U_- &= \begin{bmatrix}
u(0)& u(1)& \dots& u(T-1)
\end{bmatrix}\in \mathbb{R}^{m\times T}, \\
X_+ &= \begin{bmatrix}
x(1)& x(2)& \dots& x(T)
\end{bmatrix}\in \mathbb{R}^{n\times T},
\end{align*}
which satisfy the system dynamics (\ref{equ:sys})
\begin{equation}\label{equ:dynamics}
X_{+} = AX_{-}+ BU_-.
\end{equation}

Throughout the paper, we assume that the following block matrix of input and state data
$$
D_- = \begin{bmatrix}
U_- \\
X_-
\end{bmatrix}\in \mathbb{R}^{(m+n)\times T}
$$
has full row rank
\begin{equation}\label{equ:rank}
	\text{rank}(D_-) = m+n,
\end{equation}
i.e., the information in the data is sufficiently rich. This condition is necessary for identifying $(A,B)$ from data and for solving the data-driven LQR problem~\cite{van2020data}. As shown in \cite{de2019formulas}, it can be ensured provided that the input data $U_-$ is PE of order $n+1$. Note that the columns of $(X_-,U_-,X_+)$ are not necessarily consecutive data samples. In fact, they could be from independent or multiple averaged experiments as long as they satisfy (\ref{equ:dynamics}) and (\ref{equ:rank})~\cite{de2019formulas}.

Under the rank condition (\ref{equ:rank}), there exists a matrix $G \in \mathbb{R}^{T\times n}$ that satisfies
\begin{equation}\label{equ:relation}
\begin{bmatrix}
K \\
I_n
\end{bmatrix}=
D_-G
\end{equation}
for any given $K$. That is, $K$ can be parameterized by $K=U_-G$ where $G$ satisfies a linear constraint $X_-G= I_n$. Then, the closed-loop matrix can be expressed in a data-driven fashion as \cite{de2019formulas}
$$
A+BK=[B~~A]\begin{bmatrix}
K \\
I_n
\end{bmatrix}=(AX_-+BU_-)G=X_+G,
$$
leading to the following closed-loop system
\begin{equation}\label{equ:ddsys}
x(t+1 )=X_+Gx(t).
\end{equation}
Furthermore, the LQR problem becomes
\begin{equation}\label{prob:equi}
\begin{aligned}
&\mathop{\text {minimize}}\limits_{G}~ J(G),\\
&\text{subject to}~G\in\mathcal{S}_G :=\{G|X_-G=I_n, \rho(X_+G)<1\}.
\end{aligned}
\end{equation}
Here, $J(G)$ is the LQR cost following (\ref{equ:ddsys}) and $u(t) = U_-Gx(t)$, and $\mathcal{S}_G$ is the feasible set. In contrast to the model-based LQR, the problem (\ref{prob:equi}) is characterized by raw data matrices. Though (\ref{prob:equi}) can be reformulated as a semi-definite program (SDP) using techniques from \cite{de2019formulas,de2021low}, it is computationally challenging to solve for a large data size.

In this paper, we take an iterative PO perspective to solve (\ref{prob:equi}) viewing $G$ as the optimization matrix. We aim to design a gradient-based method to find an optimal $G$ while maintaining feasibility, and recover the control from (\ref{equ:relation}) as $K = U_-G$. Since (\ref{prob:equi}) is a challenging constrained nonconvex problem, we leverage a novel convex parameterization to establish the global convergence.



\section{Data-enabled policy optimization}\label{sec:deepo}
In this section, we first present our novel PO method for solving (\ref{prob:equi}). Then, we propose a convex parameterization of (\ref{prob:equi}) to derive the projected gradient dominance property of $J(G)$. By establishing that $J(G)$ is locally smooth over any sublevel set, we are able to show the global convergence of our method.

\subsection{Data-enabled policy optimization to solve (\ref{prob:equi})}

For $G \in \mathcal{S}_G$, the cost $J(G)$ is finite and has the following closed-form expressions \cite{de2019formulas}
\begin{equation}\label{def:JG}
J(G)= \text{Tr}\{P_G\} =\text{Tr}\{(Q+G^{\top}U_-^{\top}RU_-G)\Sigma_G\},
\end{equation}
where $P_G$ satisfies the Lyapunov equation
\begin{equation}\label{equ:Lya_P}
P_G = Q + G^{\top}U_-^{\top}RU_-G + G^{\top}X_{+}^{\top}P_GX_{+}G,
\end{equation}
and $\Sigma_G:=\mathbb{E}_{x(0)\sim \mathcal{D}}[\sum_{t=0}^{\infty}x(t)x(t)^{\top}]$ is the state covariance matrix of the closed-loop system (\ref{equ:ddsys}) satisfying
$$
\Sigma_G = I_n + X_{+}G\Sigma_G G^{\top}X_{+}^{\top}.
$$

We have the following gradient expression for $J(G)$.
\begin{lemma}\label{lem:gradient}
	For $G \in \mathcal{S}_G$, the gradient of $J(G)$ is
	$$
	\nabla J(G) = 2 E_G \Sigma_G
	$$
	with $E_G := (U_-^{\top}RU_-+X_{+}^{\top}P_GX_{+})G$.
\end{lemma}
\begin{proof}
	The proof follows from standard matrix analysis~\cite{horn2012matrix} and is similar to that of \cite[Lemma 1]{fazel2018global}.
\end{proof}

The expression of $\nabla J(G)$ is data-driven since both $E_G$ and $\Sigma_G$ can be computed using raw data matrices under the rank condition (\ref{equ:rank}).

The feasible set $\mathcal{S}_G$ contains a linear constraint $X_-G=I_n$, which motivates the use of projected gradient methods to ensure feasibility. Define the nullspace of $X_-$ as 
$$\mathcal{N}(X_-):=\{G\in \mathbb{R}^{T \times n} |X_-G = 0\},$$ 
and the projection operator $\Pi_{X_-}:= I_T-X_-^{\dagger}X_-$ onto $\mathcal{N}(X_-)$. The projected gradient update is then given by
\begin{equation}\label{equ:gd}
G^+ = G - \eta\Pi_{X_-}\nabla J(G),
\end{equation}
where $\eta \geq 0$ is the stepsize. We refer to this method as data-enabled policy optimization (DeePO) since the update (\ref{equ:gd}) can be efficiently computed by raw data matrices, and the control can be recovered from (\ref{equ:relation}) as $K = U_-G$.  As an iterative search method, the initial policy $G^0$ requires to satisfy $G^0 \in \mathcal{S}_G$.

Due to non-convexity of both the objective $J(G)$ and the constraint $\mathcal{S}_G$, it is challenging to provide global convergence guarantees for DeePO. Moreover, an optimal solution to (\ref{prob:equi}) is not unique. In fact, it has been shown in \cite[Lemma 2.1]{dorfler22on} that the solution set is 
\begin{equation}\label{equ:solset}
\left\{G|G = G^* + \Delta, \Delta \in \mathcal{N}(D_-) \right\}
~\text{with}~
G^* = D_-^{\dagger}\begin{bmatrix}
K^* \\
I_n
\end{bmatrix},
\end{equation}
which contains a considerable nullspace. Nevertheless, based on a recent work~\cite{sun2021learning} that proves optimality via convex parameterization, we are able to show a projected gradient dominance property of $J(G)$. 

\subsection{Optimality via a convex parameterization}
We first relate (\ref{prob:equi}) to a convex parameterization via a change of variables $G=L\Sigma^{-1}$ as
\begin{equation}\label{prob:convex}
\begin{aligned}
&\mathop{\text{minimize}}\limits_{L,\Sigma} ~f(L,\Sigma):= \text{Tr}\{Q\Sigma\} + \text{Tr}\{L\Sigma^{-1}L^{\top}U_-^{\top}RU_-\}, \\
&\text{subject to} ~~\Sigma = X_-L,~
\begin{bmatrix}
\Sigma-I_n & X_{+}L \\
L^{\top}X_{+}^{\top} & \Sigma
\end{bmatrix} \succeq 0.
\end{aligned}
\end{equation}	
Let $\mathcal{S}$ be its feasible set. The equivalence between the two problems (\ref{prob:equi}) and (\ref{prob:convex}) are established below.

\begin{lemma}\label{lem:asmp1}
	For any $(L,\Sigma) \in \mathcal{S}$, $\Sigma$ is invertible and $L\Sigma^{-1} \in \mathcal{S}_G$. Moreover, for $G\in \mathcal{S}_G$ it holds that
	\begin{equation}\label{equ:equi}
	J(G) = \min_{L,\Sigma} \{f(L,\Sigma),\text{s.t.} (L,\Sigma) \in \mathcal{S}, L\Sigma^{-1} = G\}.
	\end{equation}
\end{lemma}
\begin{proof}
	Applying the Schur complement to the LMI constraint in \eqref{prob:convex} yields $\Sigma \succ 0$ and
	$$
	\Sigma-I_n -X_{+}L\Sigma^{-1}L^{\top}X_{+}^{\top} \succeq 0.
	$$
	Due to non-singularity of $\Sigma$, let $G = L\Sigma^{-1}$. Then, a substitution of $L = G\Sigma$ into the above inequality yields
	$$
	\Sigma-I_n -X_{+}G\Sigma G^{\top}X_{+}^{\top} \succeq 0.
	$$
	Thus, $X_{+}G$ is stable, i.e., $\rho(X_{+}G)<1$. Since the first constraint of \eqref{prob:convex} implies $X_-G = X_-L\Sigma^{-1} = \Sigma \Sigma^{-1} = I_n$, it holds that $G = L\Sigma^{-1} \in \mathcal{S}_G$.
	
	Next, we prove the second statement. Using the constraint $G = L\Sigma^{-1}$ and the Schur complement, the right-hand side of (\ref{equ:equi}) becomes
	\begin{equation}\label{prob:JG}
	\begin{aligned}
	&\min_{\Sigma} ~ \text{Tr}((Q+G^{\top}U_-^{\top}RU_-G)\Sigma)\\
	&~\text{s.t.} ~X_-G = I_n, \Sigma \succ 0, \Sigma \succeq I_n + X_{+}G\Sigma G^{\top}X_{+}^{\top}.
	\end{aligned}
	\end{equation}
	Let $\Sigma(\Theta)$ be the unique positive definite solution of the Lyapunov equation
	$$
	\Sigma(\Theta) = \Theta + X_{+}G\Sigma(\Theta) G^{\top}X_{+}^{\top}
	$$
	with $\Theta \succeq I_n$. By monotonicity of $\Sigma(\Theta)$, we have $\Sigma(\Theta) \succeq \Sigma(I_n)$. Since $Q+G^{\top}U_-^{\top}RU_-G \succ 0$, the minimum of (\ref{prob:JG}) is attained at $\Sigma(I_n)$, which is
	$
	\text{Tr}((Q+G^{\top}U_-^{\top}RU_-G)\Sigma(I_n))
	$
	with $X_-G = I_n$. This is the definition of $J(G)$ in (\ref{def:JG}). 
\end{proof}

In the following lemma, we show the convexity of the parameterization \eqref{prob:convex}.
{\begin{lemma}\label{lem:asmp2}
	The feasible set $\mathcal{S}$ of (\ref{prob:convex}) is convex in $(L,\Sigma)$, and $f(L,\Sigma)$ is differentiable over an open domain that contains $\mathcal{S}$. Moreover, $f(L,\Sigma)$ is convex over $\mathcal{S}$.
\end{lemma}}
\begin{proof}
	Since the constraints in (\ref{prob:convex}) are linear in $(L,\Sigma)$, the feasible set $\mathcal{S}$ is convex. Clearly,  $f(L,\Sigma)$ is differential over $\mathcal{S}$. Define the Hessian operator acting on the direction $(\tilde{L},\tilde{\Sigma})$:
	$$
	h(L,\Sigma;\tilde{L},\tilde{\Sigma}):= \nabla^2 f(L,\Sigma)[(\tilde{L},\tilde{\Sigma}),(\tilde{L},\tilde{\Sigma})],
	$$
	which by standard matrix analysis~\cite{horn2012matrix} can be written as
	$$
	h(L,\Sigma;\tilde{L},\tilde{\Sigma})=2\left\|R^{\frac{1}{2}}(U_-\tilde{L}- U_-L\Sigma^{-1} \tilde{\Sigma}) \Sigma^{-\frac{1}{2}}\right\|_F^2\geq 0.
	$$	
	Thus, $f$ is convex over $\mathcal{S}$.
\end{proof}

We now formally define the gradient dominance property.{\begin{definition} A differentiable function $g(x): \mathbb{R}^n \rightarrow \mathbb{R}$ with a finite global minimum $g^*$ is gradient dominated of degree $p$ over a set $\mathcal{X}\subseteq \text{dom}(g) $ if
$$
g(x) - g^* \leq \lambda_{\mathcal{X}} \| \nabla g(x) \|^p, ~~\forall x \in \mathcal{X},~\text{for some}~ \lambda_\mathcal{X} > 0.
$$
\end{definition}} 

The gradient dominance property means that all the stationary points are optimal. Moreover, the convergence rate of gradient-based methods usually depends on the values of the degree $p$. Particularly, for smooth objective function $p=1$ leads to a sublinear rate and $p=2$ leads to a linear rate.

Equipped with Lemmas \ref{lem:asmp1} and \ref{lem:asmp2}, we apply \cite[Theorem 1]{sun2021learning} to show the gradient dominance property of $J(G)$ over any sublevel set
$
S_G(a):= \{G\in \mathbb{R}^{T\times n} | J(G)\leq a\}
$
with $a > 0$.
{\begin{lemma}[\textbf{Projected gradient dominance of degree 1}]\label{lem:pl}
	For $G \in \mathcal{S}_G(a)$, there exists $\mu(a) >0$ such that
	$$
	J(G) - J^* \leq \mu(a) \|\Pi_{X_-} \nabla J(G)\|,
	$$
	where $J^*$ is the optimal LQR cost to (\ref{prob:equi}).
\end{lemma}}
\begin{proof}
	By Lemmas \ref{lem:asmp1} and \ref{lem:asmp2}, the data-driven LQR problem \eqref{prob:equi} and its convex parameterization \eqref{prob:convex} satisfy the assumptions required to apply \cite[Theorem 1]{sun2021learning}. Then, there exists $c(a)>0$ and a direction $V\in \mathcal{N}(X_-)$ with $\|V\|_F=1$ in the descent cone of $\mathcal{S}_G(a)$ such that
	$$
	J'(G)[V] \leq -c(a) (J(G)-J^*),
	$$
	where $J'(G)[V]$ denotes the derivative along the direction $V$. Let 
	$V' = \Pi_{X_-} \nabla J(G) / \|\Pi_{X_-} \nabla J(G)\|_F$ be the normalized projected gradient. Then, we have $ J'(G)[V'] \leq  J'(G)[V]$ since both $V$ and $V'$ are in $\mathcal{N}(X_-)$, and $V'$ is the direction of the projection of the gradient. Thus, we have $J(G) - J^* \leq \mu(a)\|\Pi_{X_-} \nabla J(G)\|$ with $\mu(a) = 1/c(a)$.
	
	Next, we derive an explicit upper bound of $\mu(a)$ over $G \in \mathcal{S}_G(a)$. By \cite[Theorem 1]{sun2021analysis}, $c(a)$ is given by
	$$
	c(a) = (2\max\{\|L-L^*\|_F/\underline{\sigma}(\Sigma), \|\Sigma - \Sigma^*\|_F  \|L\| /\underline{\sigma}^{2}(\Sigma) \})^{-1},
	$$
	where $(L^*,\Sigma^*)$ is an optimal point and $(L,\Sigma) = \arg\min_{L',\Sigma'} f(L',\Sigma')$ subject to $(L',\Sigma') \in \mathcal{S}, L'(\Sigma')^{-1} = G$. We now provide upper bounds for $\underline{\sigma}^{-1}(\Sigma), \|L\|_F, \|\Sigma\|_F$. Since $\underline{\sigma}(\Sigma)\geq 1$, it holds $\underline{\sigma}^{-1}(\Sigma) \leq 1$. The sublevel set gives
	$\text{Tr}\{Q\Sigma\} \leq a$, and hence $\|\Sigma\|_F \leq a / \underline{\sigma}(Q)$. Since
	\begin{align*}
		&\underline{\sigma}(R) \underline{\sigma}^2({U}_-)\|\Sigma\|^{-1} \|L\|^2_F \leq \text{Tr}\{L\Sigma^{-1}L^{\top}{U}_-^{\top}R{U}_-\}  \\
		&\leq \text{Tr}\{Q\Sigma\} + \text{Tr}\{L\Sigma^{-1}L^{\top}{U}_-^{\top}R{U}_-\} \leq a,
	\end{align*}
	an upper bound of $\|L\|_F$ is given by
	$$
	\|L\|_F \leq  \left(\frac{a\|\Sigma\|}{\underline{\sigma}(R) \underline{\sigma}^2({U}_-)}\right)^{{1}/{2}} \leq \frac{a}{(\underline{\sigma}(Q)\underline{\sigma}(R))^{1/2} \underline{\sigma}({U}_-)}.
	$$
	Those bounds are also true for $L^*,\Sigma^*$. Furthermore, we can provide an upper bound of $\mu(a)$ as 
	$$
	\mu(a) \leq \frac{4a}{(\underline{\sigma}(Q)\underline{\sigma}(R))^{1/2} \underline{\sigma}({U}_-)} \max\left\{1, \frac{a}{\underline{\sigma}(Q)} \right\}.
	$$
	The proof is completed.
\end{proof}

In contrast to the existing literature~\cite{fazel2018global} on PO for the LQR, the cost $J(G)$ here is \textit{projected gradient} dominated, meaning that $G$ is optimal if the projected gradient $\Pi_{X_-} \nabla J(G)$ is equal to zero. By using Lemma \ref{lem:pl}, we next show global convergence of the projected gradient descent in (\ref{equ:gd}). 

\subsection{Global convergence of DeePO}
We first prove the smoothness of $J(G)$. Since $J(G)$ tends extremely to infinity as $G$ approaches the boundary $\partial \mathcal{S}_G$, we can only show that $J(G)$ is \textit{locally} smooth over any sublevel set. Define the Hessian acting on the direction $Z\in \mathbb{R}^{T\times n}$ as
$
\nabla^2 J(G)[Z,Z] := \left.\frac{d^2}{dt^2}J(G+tZ)\right|_{t=0},
$
and the directional derivative of $P_G$ as
$
P_G'[Z]:=\left.\frac{d}{dt}P_{G+tZ}\right|_{t=0}.
$ Then, we have the following closed-form expression for the Hessian.
\begin{lemma}
	For $G \in \mathcal{S}_G$ and a feasible direction $Z\in \mathbb{R}^{T\times n}$, the Hessian of  $J(G)$ is characterized by 
	\begin{align*}
	\nabla^2 J(G)[Z,Z]&=2\text{Tr}\{Z^{\top}(U_-^{\top}RU_-+X_{+}^{\top}P_GX_{+})Z\Sigma_G  \}\\
	&~~~+ 4\text{Tr}\{Z^{\top}X_+^{\top}P_G'[Z]X_{+}G\Sigma_G  \},
	\end{align*}
	where $P_G'[Z] = \sum_{i=0}^{\infty} (G^{\top}X_{+}^{\top})^i(Z^{\top}E_G + E_G^{\top}Z)(X_{+}G)^i$.
\end{lemma}
\begin{proof}
	The proof follows from standard matrix analysis~\cite{horn2012matrix} and is omitted due to space limitation.
\end{proof}

Define $\|\nabla^2 J(G)\| := \sup_{\|Z\|_F=1} \left|\nabla^2 J(G)[Z,Z]\right|$. We show an upper bound for $\|\nabla^2 J(G)\|$ over a sublevel set.

\begin{lemma}[\textbf{Local smoothness}]\label{lem:smooth}
	For $G \in \mathcal{S}_G(a)$, it holds
	$$
	\|\nabla^2 J(G)\| \leq \text{poly}(a,\|U_-\|,\|X_+\|_F,\|R\|,\underline{\sigma}(Q)):=l(a),
	$$
	where $l(a)$ is the smoothness constant of $J(G)$ over $\mathcal{S}_G(a)$. That is, for any $G,G' \in \mathcal{S}_G(a)$ satisfying $G+\delta(G'-G) \in \mathcal{S}_G(a), \forall \delta \in [0,1]$, the following inequality holds
	$$
	J(G') \leq J(G) +  \langle \nabla J(G), G' - G \rangle + {l(a)}\|G'-G\|^2/2.
	$$
\end{lemma}

The proof is technical and provided in Appendix \ref{app:2}. 

Under the gradient dominance property of degree 1 in Lemma \ref{lem:pl} and the local smoothness in Lemma \ref{lem:smooth}, we now show the global sublinear convergence of DeePO. The key is to select an appropriate stepsize such that the policy sequence is feasible and stays in the sublevel set associated with the initial policy $G^0 \in \mathcal{S}_G$. For simplicity, let $\mu_0$ and $l_0$ denote the projected gradient dominance and smoothness constants of $J(G)$ over $\mathcal{S}_G(J(G^0))$, respectively. We present our convergence result in the following theorem.
{
\begin{theorem}[\textbf{Global convergence}]\label{thm:conv}
	For $G^0\in \mathcal{S}_G$ and a stepsize $\eta \in (0,1/l_0]$, the update (\ref{equ:gd}) leads to $G^k \in \mathcal{S}_G(J(G^0)), \forall k \in \mathbb{N}$. Moreover, for any $\epsilon> 0$ and
	\begin{equation}\label{equ:kk}
	k \geq \frac{2\mu_0^2}{\epsilon(2\eta - l_0\eta^2)},
	\end{equation}
	the update (\ref{equ:gd}) enjoys the following performance bound
	$$
	J(G^k)-J^* \leq \epsilon.
	$$
\end{theorem}}
\begin{proof}
	Define $G_\eta := G- \eta \Pi_{X_-}\nabla J(G)$. We first show that for a non-optimal $G\in\mathcal{S}_G(a)$ and any $\eta \in [0,1/l(a)]$, it holds $G_{\eta} \in \mathcal{S}_G(a)$. 
	
	Define $\mathcal{S}_G^o(a):= \{G\in \mathcal{S}_G|J(G)<a\}$, and its complement as $(\mathcal{S}_G^{o}(a))^c$, which is closed. By Lemma \ref{lem:smooth}, given $\phi\in (0,1)$, there exists $b>0$ such that $\|\nabla^2 J(G)\|\leq (1+\phi)l(a)$ for $G \in \mathcal{S}_G(a+b)$. Clearly, $\mathcal{S}_G(a) \cap (\mathcal{S}_G^{o}(a+b))^c = \emptyset$. Then, the distance between them $d:=\inf\{\|G'-G\|, \forall G \in \mathcal{S}_G(a), G'\in (\mathcal{S}_G^{o}(a+b))^c \}$ is positive. 
	
	Let $\overline{N} \in \mathbb{N}_+$ be large enough such that $2/(\overline{N}(1+\phi)l(a)) < d/\|\Pi_{X_-}\nabla J(G)\|$, which is well-defined since $G$ is not optimal. Define a stepsize $\tau \in [0, 2/(\overline{N}(1+\phi)l(a))]$. Since $\tau < d/\|\Pi_{X_-}\nabla J(G)\|$, we have $\|G_{\tau}-G\|< d$, i.e., $G_{\tau}\in \mathcal{S}_G(a+b)$. Thus, we can apply Lemma \ref{lem:smooth} over $\mathcal{S}_G(a+b)$ to show
	$$
		J(G_{\tau}) -  J(G) \leq - \tau(1 -\frac{(1+\phi)l(a)\tau}{2})\| \Pi_{X_-} \nabla J(G)\|^2\leq 0,
	$$
	where the last inequality follows from $\tau \leq 2/(1+\phi)l(a)$. This implies that the segment between $G$ and $G_{\tau}$ is contained in $\mathcal{S}_G(a)$. It is also clear that $G_{2\tau}\in \mathcal{S}_G(a+b)$ since $\|G_{2\tau}-G_{\tau}\|< d$. Then, we can use induction to show that the segment between $G$ and $G_{N\tau}$ for $N\in \mathbb{N}_+$ is in $ \mathcal{S}_G(a)$ as long as $N\tau \leq 2/(1+\phi)l(a)$. Since $\phi \in (0,1)$, we let $\eta \leq 1/l(a)$ to ensure the segment between $G$ and $G_{\eta}$ to be contained in $ \mathcal{S}_G(a)$. 

	Then, a simple induction leads to that for $\eta \in [0,1/l_0]$, the update (\ref{equ:gd}) satisfies $G^k \in \mathcal{S}_G(J(G^0)), \forall k \in \mathbb{N}$. Moreover, the cost satisfies
	$$
	J(G^{k+1}) \leq  J(G^k) - \eta(1 -\frac{l_0\eta}{2})\| \Pi_{X_-} \nabla J(G^k)\|^2.
	$$	
	Using Lemma \ref{lem:pl} and subtracting $J^*$ in both sides yields
	$$
	J(G^{k+1})-J^* \leq J(G^k) - J^* - \frac{2\eta - l_0\eta^2}{2\mu_0^2}(J(G^k) - J^*)^2
	$$ 
	Let $e^k = J(G^k) - J^*$. Dividing by $e^ke^{k+1}$ in both sides and noting $e^{k+1}\leq e^k$ leads to
	$$
	\frac{2\eta - l_0\eta^2}{2\mu_0^2} \leq \frac{1}{e^{k+1}} - \frac{1}{e^k}
	$$
	Summing up both sides over $0,1,\dots,k-1$ and using telescopic cancellation yields that
	$$
	\frac{k(2\eta - l_0\eta^2)}{2\mu_0^2}\leq \frac{1}{e^k} - \frac{1}{e^0} \leq \frac{1}{e^k}
	$$ 
	Letting the right-hand side of the above inequality equal $\epsilon$ and solving $k$ yields (\ref{equ:kk}) under $\eta \in (0,1/l_0]$.
\end{proof}

We compare with the traditional PO for the LQR \cite{fazel2018global,mohammadi2022convergence,malik2019derivative}. Their approach relies on a zeroth-order estimate of the policy gradient, which inevitably requires numerous system trajectories to approximate the cost. In sharp contrast, DeePO directly computes the gradient from a batch of raw data matrices based on a data-based representation of the closed-loop system. This remarkable feature enables DeePO to work with only a small set of PE data. Moreover, the state-of-the-art sample complexity (in terms of number of sampled trajectories, the length of which can be very long) of PO in \cite{fazel2018global,mohammadi2022convergence,malik2019derivative} is $\mathcal{O}(\log(1/\epsilon))$, while our sample complexity (in terms of number of state-input pairs) is independent of $\epsilon$. Even though both two approaches achieve global convergence (albeit with vastly different amounts of data), DeePO is more flexible as it is compatible with regularization methods used to enhance the robustness to noisy data, which will be shown in the next section. To the best of our knowledge, there are no robustifying regularization methods that have been applied to the PO method for the LQR problem.

\section{DeePO for the regularized LQR}\label{sec:reg}
For the direct data-driven LQR formulation~\cite{de2021low,dorfler2021certainty,dorfler22on}, regularization plays an important role in promoting certainty-equivalence and robust stability when the data is corrupted with noise. This section investigates how regularization affects the convergence of DeePO.

\subsection{Certainty-equivalence regularizer}
Consider the regularized LQR problem
\begin{equation}\label{prob:cq}
\begin{aligned}
&\mathop{\text {minimize}}\limits_{G}~J_{\lambda}(G) := J(G)+\lambda \|\Pi_{D_-}G\Sigma_G^{1/2}\|^2,~~\\
&\text{subject to} ~~ G \in \mathcal{S}_G,
\end{aligned}
\end{equation}
where $\lambda\geq0$ is a user-defined constant and $\Pi_{D_-} := I - D_-^{\dagger}D_-$ is the projection matrix onto the nullspace of $D_-$. For the noiseless data $(X_-,U_-,X_+)$ here, the orthogonality regularizer in (\ref{prob:cq}) does not change the optimal cost but only singles out a solution $G^*$ satisfying $\Pi_{D_-}G^* = 0$ from the solution set in (\ref{equ:solset}). When the data is corrupted with noises, it promotes certainty-equivalence, i.e., when $\lambda$ tends to infinity the solution of (\ref{prob:cq}) coincides with that of indirect data-driven control with an underlying maximum likelihood system identification attenuating the effect of noise; we refer interested readers to \cite[Section III]{dorfler22on} for more discussions. 

Note that we have added the weighting $\Sigma_G^{1/2}$ to the regularizer (c.f. \cite[(15)]{dorfler22on}) to make it compatible with the convex parameterization (\ref{prob:convex}). As a result, (\ref{prob:cq}) can be formulated with $L\Sigma^{-1}=G$ as the following convex problem
\begin{equation}\label{prob:reg1}
\begin{aligned}
&\mathop{\text {minimize}}\limits_{L,\Sigma} ~f_{\lambda}(L,\Sigma):= \text{Tr}\{Q\Sigma\} \\
&~~~~~~~~~~~~+\text{Tr}\{L\Sigma^{-1}L^{\top}(\lambda\Pi_{D_-}^{\top}\Pi_{D_-}+U_-^{\top}RU_-)\},\\
&\text{subject to} ~~\Sigma = X_-L,~
\begin{bmatrix}
\Sigma-I_n & X_{+}L \\
L^{\top}X_{+}^{\top} & \Sigma
\end{bmatrix} \succeq 0.
\end{aligned}
\end{equation}


Comparing (\ref{prob:reg1}) with (\ref{prob:convex}), we see that $f_{\lambda}(L,\Sigma)$ upon amounts to $f(L,\Sigma)$ adding a convex regularizer, and hence $f_{\lambda}(L,\Sigma)$ is convex. Indeed, by standard matrix analysis~\cite{horn2012matrix}, its Hessian acting on the direction $(\tilde{L},\tilde{\Sigma})$ satisfies
\begin{align*}
&\nabla^2 f_\lambda(L,\Sigma)[(\tilde{L},\tilde{\Sigma}),(\tilde{L},\tilde{\Sigma})] =\nabla^2 f(L,\Sigma)[(\tilde{L},\tilde{\Sigma}),(\tilde{L},\tilde{\Sigma})]\\
&+	2\lambda\|(\Pi_{D_-}\tilde{L}- \Pi_{D_-}L\Sigma^{-1} \tilde{\Sigma}) \Sigma^{-1/2}\|_F^2\\
&\geq \nabla^2 f(L,\Sigma)[(\tilde{L},\tilde{\Sigma}),(\tilde{L},\tilde{\Sigma})].
\end{align*}
Moreover, following analogous arguments as in Section \ref{sec:deepo}, $J_{\lambda}(G)$ can also be shown to be locally smooth. Based on previous analysis, the projected gradient update
\begin{equation}\label{equ:pgregu}
G^+ = G - \eta\Pi_{X_-}\nabla J_{\lambda}(G)
\end{equation}
converges to the optimal solution of (\ref{prob:cq}) under a proper stepsize selection.

\subsection{Robustness-promoting regularizer}\label{subsec:rb}
Regularization can also be used to enhance robust stability. Consider the following regularized LQR problem
\begin{equation}\label{prob:rob}
\begin{aligned}
&\mathop{\text {minimize}}\limits_{G}~ J_{\gamma}(G) :=  J(G)+\gamma \text{Tr}\{G\Sigma_GG\},~~\\
&\text{subject to} ~ G \in \mathcal{S}_G,
\end{aligned}
\end{equation}
where $\gamma\geq 0$ is a user-defined constant. To see why it promotes the robust stability for noisy data, we note that the state covariance matrix is given by
$$
\Sigma_G = I_n + X_{+}G\Sigma_G G^{\top}X_{+}^{\top}.
$$
Thus, a small $\text{Tr}\{G\Sigma_GG^{\top}\}$ can reduce the effect of noises in $X_+$. Different from the certainty-equivalence regularization, the regularizer in (\ref{prob:rob}) bias the LQR solution even when the data is noiseless, reflecting a trade-off between performance and robustness.

The problem (\ref{prob:rob}) can be formulated with $L\Sigma^{-1}=G$ as
\begin{equation}\label{prob:reg2}
\begin{aligned}
&\mathop{\text {minimize}}\limits_{L,\Sigma} ~f_{\gamma}(L,\Sigma):=  \text{Tr}\{Q\Sigma\}\\
&~~~~~~~~~~~~~+\text{Tr}\{L\Sigma^{-1}L^{\top}(\gamma I_T+U_-^{\top}RU_-)\},\\
&\text{subject to} ~~\Sigma = X_-L,~
\begin{bmatrix}
\Sigma-I_n & X_{t+1}L \\
L^{\top}X_{t+1}^{\top} & \Sigma
\end{bmatrix} \succeq 0.
\end{aligned}
\end{equation}	

Clearly, $f_{\lambda}(L,\Sigma)$ is also convex since
\begin{align*}
&\nabla^2 f_\gamma(L,\Sigma)[(\tilde{L},\tilde{\Sigma}),(\tilde{L},\tilde{\Sigma})] =\nabla^2 f(L,\Sigma)[(\tilde{L},\tilde{\Sigma}),(\tilde{L},\tilde{\Sigma})]\\
&+	2\gamma\|(\tilde{L}- L\Sigma^{-1} \tilde{\Sigma}) \Sigma^{-\frac{1}{2}}\|_F^2
\geq \nabla^2 f(L,\Sigma)[(\tilde{L},\tilde{\Sigma}),(\tilde{L},\tilde{\Sigma})].
\end{align*}
By analogous reasoning and combining the smoothness of the regularizer, the projected gradient update
\begin{equation}\label{equ:rb}
	G^+ = G - \eta\Pi_{X_-}\nabla J_{\gamma}(G)
\end{equation}
converges to the optimal solution of (\ref{prob:rob}) under a proper stepsize selection.

\subsection{Implicit regularization}

Apart from the convergence, we observe an interesting \textit{implicit regularization} property of the certainty-equivalence regularized LQR problem (\ref{prob:cq}) formally defined below.

\begin{definition}[\textbf{Implicit regularization}]
	For the regularized LQR problem (\ref{prob:cq}), suppose that a convergent algorithm generates a sequence of $\{G^k\}$. If $G^{\infty}:=\lim\limits_{k\rightarrow\infty}G^k$ satisfies $\Pi_{D_-}G^{\infty}=0$, then the algorithm is called \textit{regularized}; If it is regularized with $\lambda=0$, then it is called \textit{implicitly} regularized.
\end{definition}

The concept of implicit regularization has been adopted in many recent works on nonconvex optimization, including deep learning~\cite{neyshabur2017geometry}, matrix factorization~\cite{arora2019implicit}, and also PO for robust LQR problems~\cite{zhang2021policy}. As its name suggests, it means that the algorithm without regularization behaves as if it is regularized. Note that implicit regularization is a property of a certain algorithm for solving a certain nonconvex problem. In the following theorem, we specify the conditions for the update (\ref{equ:pgregu}) to be implicitly regularized for problem (\ref{prob:cq}).
\begin{figure}[t]
	\centerline{\includegraphics[width=50mm]{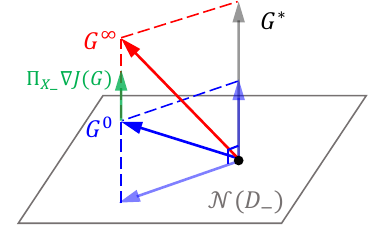}}
	\caption{Supspace relations among $\mathcal{N}(D_-)$, $\Pi_{X_-} \nabla J(G)$, and $G^*$.}
	\label{pic:subspace}
\end{figure}

\begin{theorem}[\textbf{Implicit regularization}]\label{thm:impli}
	Consider (\ref{prob:cq}) with $\lambda = 0$ and suppose that $G^0$ satisfies $\Pi_{D_-}G^0 = 0$. Then, the update (\ref{equ:pgregu}) leads to $\Pi_{D_-}G^k = 0, k\in \{0,1,\dots\}$.
\end{theorem}
\begin{proof}
	Since $\lambda = 0$, it suffices to show that $\Pi_{X_-}\nabla J(G)$ is orthogonal to the nullspace of $D_-$.
	
	By using Lemma \ref{lem:gradient}, the gradient of $J(G)$ is written as
	\begin{align*}
	\nabla J(G) & = 2(U_-^{\top}RU_-+X_{+}^{\top}P_GX_{+})G \Sigma_G \\
	& = 2 \begin{bmatrix}
	U_- \\
	X_-
	\end{bmatrix}^{\top}\begin{bmatrix}
	R+B^{\top}PB & B^{\top}PA \\
	A^{\top}PB & A^{\top}PA
	\end{bmatrix}\begin{bmatrix}
	U_- \\
	X_-
	\end{bmatrix}G\Sigma_G.
	\end{align*}
	
	
	We also have the following observation
	\begin{align*}
	&(I - X_-^{\dagger}X_-)\begin{bmatrix}
	U_- \\
	X_-
	\end{bmatrix}^{\top}\\
	 &= \begin{bmatrix}
	U_-^{\top}-X_-^{\dagger}X_-U_-^{\top} & 0
	\end{bmatrix} \\
	& = \begin{bmatrix}
	U_-^{\top}-X_-^{\top}(X_-X_-^{\top})^{-1}X_-U_-^{\top} & 0 
	\end{bmatrix} \\
	& = \begin{bmatrix}
	U_- \\
	X_-
	\end{bmatrix}^{\top} \begin{bmatrix}
	I_m & 0 \\
	-(X_-X_-^{\top})^{-1}X_-U_-^{\top} & 0
	\end{bmatrix}.
	\end{align*}
	
	Thus, $\Pi_{X_-}\nabla J(G)$ is in the range space of $D_- = \begin{bmatrix}
	U_-^{\top} &
	X_-^{\top}
	\end{bmatrix}^{\top}$, and hence
	$
	\Pi_{D_-} \Pi_{X_-} \nabla J(G) = 0.
	$ The update (\ref{equ:pgregu}) further leads to $\Pi_{D_-} G^{k+1} = \Pi_{D_-} G^k - \eta \Pi_{D_-} \Pi_{X_-} \nabla J(G^k) = \Pi_{D_-} G^k=0$.
\end{proof}

By Theorem \ref{thm:impli}, a sufficient condition for implicit regularization is 
$$
G^0 = D_-^{\dagger}\begin{bmatrix}
K^0 \\
I_n
\end{bmatrix},
$$
provided with a stabilizing policy $K^0$. Theorem \ref{thm:impli} also helps understand the optimization landscape of DeePO. Fig. \ref{pic:subspace} illustrates the relations among the nullspace $\mathcal{N}(D_-)$, the projected gradient, and an optimal solution $G^*$. Since $\Pi_{X_-} \nabla J(G)$ is orthogonal to $\mathcal{N}(D_-)$, the resulted policy of DeePO can be read as
$
G^{\infty} = \Pi_{D_-}G^0 + G^*.
$

\section{Simulations}\label{sec:simu}
In this section, we perform simulations to validate the convergence of DeePO and the effects of regularization.
\subsection{Numerical example}
We randomly generate a dynamical model $(A,B)$ with $n=4,m=2$ from a standard normal distribution and normalize $A$ such that $\rho(A)=0.8$, i.e., the open-loop system is stable. The resulting model parameters $(A,B)$ are
\begin{align*}
&A = \begin{bmatrix}
-0.137  &  0.146  &  -0.297  & 0.283\\
0.487  & 0.095 &  0.417 &  0.301\\
-0.018  &  0.049  &  0.175  & 0.435\\
0.143  &  0.317 &  -0.293 &  -0.107
\end{bmatrix},\\
&B = \begin{bmatrix}
1.639  &  0.930\\
0.264 &  1.793\\
-1.464  & -1.183\\
-0.776  & -0.111
\end{bmatrix}.
\end{align*}
It is straightforward to check that $(A,B)$ is controllable. Let $Q = I_4$ and $R = I_2$. We use Gaussian distribution to generate a batch of sufficiently exciting data $(U_-,X_-)$ with $T=10$ that satisfies (\ref{equ:rank}), and compute $X_+$ by (\ref{equ:dynamics}). In the sequel, we only use $(U_-,X_-,X_+)$ to perform the DeePO methods and validate the convergence.

\begin{figure}[t]
	\centerline{\includegraphics[width=60mm]{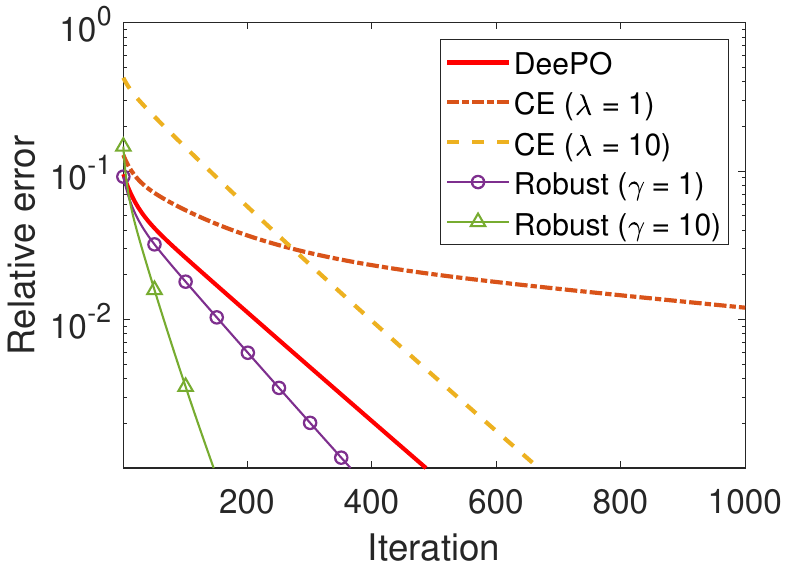}}
	\caption{Convergence of the DeePO methods.}
	\label{pic:conv}
\end{figure}
\subsection{Convergence of the DeePO methods}

We consider three algorithms, i.e, DeePO in (\ref{equ:gd}), DeePO with the certainty-equivalence regularizer in (\ref{equ:pgregu}) and with the robustness regularizer in (\ref{equ:rb}). For all the three algorithms, we set the stepsize to $\eta = 2\times 10^{-3}$ for a fair comparison. For DeePO and DeePO with robustness regularizer, we set the initial policy as
$$
G^0 = D_-^{\dagger}\begin{bmatrix}
K^0 \\ 
I_4
\end{bmatrix}
\in \mathcal{S}_G
$$
with $K^0=0$ since the system is open-loop stable. For DeePO with certainty-equivalence regularizer, we set 
$$
G^0 = D_-^{\dagger}\begin{bmatrix}
0 \\ 
I_4
\end{bmatrix} +  \Pi_{D_-}M \in \mathcal{S}_G,
$$
where the elements of $M \in \mathbb{R}^{T \times n}$ are randomly sampled from a Gaussian distribution $\mathcal{N}(0, 0.01)$ (otherwise due to the implicit regularization, there will be no difference in the convergence curve compared with DeePO). To see how regularization parameters affect the performance, we select $\lambda = 1,10$ for the certainty-equivalence regularizer and $\gamma = 1,10$ for the robustness regularizer.

We illustrate the performance of the three algorithms in Fig. \ref{pic:conv}, where their relative errors are defined as $(J(G^k)-J^*)/J^*$, $(J_{\lambda}(G^k)-J_{\lambda}^*)/J_{\lambda}^*$, and $(J_{\gamma}(G^k)-J_{\gamma}^*)/J_{\gamma}^*$, respectively. {While Theorem \ref{thm:conv} only shows a more conservative sublinear convergence rate, all the three algorithms converge linearly in the simulation.} The DeePO algorithm with certainty-equivalence regularizer (denoted by CE in Fig. \ref{pic:conv}) has the slowest convergence. The case for $\lambda = 10$ converges faster than the case $\lambda = 1$ due to the faster decay of the regularizer $\lambda \|\Pi_{D_-}G\Sigma_G^{1/2}\|^2$, and it achieves the same rate as the unregularized DeePO algorithm. Under the robustness regularizer, the DeePO algorithm has the fastest convergence, and $\gamma =10$ leads to a larger convergence rate. Nevertheless, the resulted policy is different from those of the other two algorithms as discussed in Section \ref{subsec:rb}. Finally, we note that all the algorithms only use $10$ pairs of state-input data to achieve an arbitrary relative error. In sharp contrast, the zeroth-order optimization method in \cite{malik2019derivative} uses $10^5$ trajectories (of manually tuned length to approximate the cost well) to achieve $0.01$ relative error for an LTI system with $m=n=3$.

\section{Conclusion}\label{sec:conc}
In this paper, we have proposed the DeePO method that only requires a finite number of PE data to solve the LQR problem. By relating the nonconvex optimization problem to a convex program, we have shown the global convergence of DeePO. Furthermore, we have shown that the regularization method can be applied to enhance certainty-equivalence and robust stability without affecting its convergence. The implicit regularization property has also provided an insightful understanding on the optimization landscape of DeePO. 

In future, it would be valuable to discover a strongly convex reparameterization of (\ref{prob:equi}), which may improve the sublinear convergence rate to linear. It would also be interesting to study DeePO in a more general setting, e.g., the LQR with noisy inputs. Since DeePO is an efficient iterative method, it is expected to be able to applied to online control, where the control performance is constantly improved by collecting more real-time data. We are also hopeful that it can be used to solve the adaptive LQR for time-varying systems.

\bibliographystyle{IEEEtran}
\bibliography{mybibfile}

\begin{appendices}
	\section{Proof of Lemma \ref{lem:smooth}}\label{app:2}
	We begin with a technical lemma.
	\begin{lemma}
		For $G \in \mathcal{S}_G$, it follows that
		$$
		\|\Sigma_G\| \leq \text{Tr}\{\Sigma_G\} \leq {J(G)}/{\underline{\sigma}(Q)}, \|P_G\| \leq J(G).
		$$
	\end{lemma}
	
	This lemma follows directly from the definition of $J(G)$ and is consistent with \cite[Lemma 13]{fazel2018global}.
	
	Let $Z \in \mathbb{R}^{T\times n}$ be a feasible direction with $\|Z\|_F=1$. Then, it follows that
	\begin{align*}
	\|\nabla^2 J(G)\| &\leq 2\|Z^{\top}(U_-^{\top}RU_-+X_{+}^{\top}P_GX_{+})Z\|\cdot \text{Tr}\{\Sigma_G\} \\
	&~~~+ 4\left|\text{Tr}\{Z^{\top}X_{+}^{\top}P_G'[Z]X_{+}G\Sigma_G \}\right|.
	\end{align*}
	
	The first term can be upper bounded by
	\begin{align*}
	&\|Z^{\top}(U_-^{\top}RU_-+X_{+}^{\top}P_GX_{+})Z\|\cdot \text{Tr}\{\Sigma_G\}\\
	& \leq (\|U_-\|^2\|R\|+\|X_{+}\|^2J(G)) \cdot \frac{J(G)}{\underline{\sigma}(Q)}.
	\end{align*}
	
	For the second term, we have that
	\begin{align*}
	&\left|\text{Tr}\{Z^{\top}X_{+}^{\top}P_G'[Z]X_{+}G\Sigma_G\}\right|\\
	& \leq \sup_{\|Z\|_F=1} \|Z^{\top}X_{+}^{\top}P_G'[Z]X_{+}G\Sigma_G^{1/2}\|_F\|\Sigma_G^{1/2}\|_F\\
	& \leq  \|X_{+}\|_F^2  \|X_{+}G\Sigma_G^{1/2}\|_F\|\Sigma_G^{1/2}\|_F \sup_{\|Z\|_F=1}\|P_G'[Z]\|_F \\
	&\leq \|X_{+}\|_F^2 \cdot \frac{J(G)}{\underline{\sigma}(Q)}\sup_{\|Z\|_F=1}\|P_G'[Z]\|_F,
	\end{align*}
	where the last inequality follows from the definition of $\Sigma_G$. 
	
	
	Thus, it suffices to bound $\|P_G'[Z]\|_F$. We have that
	\begin{align*}
	&Z^{\top}E_G + E_G^{\top}Z \\
	&\preceq Z^{\top}(U_-^{\top}RU_-+X_{+}^{\top}P_GX_{+})Z \\
	&~~~+ G^{\top}(U_-^{\top}RU_-+X_{+}^{\top}P_GX_{+})G\\
	&=Z^{\top}(U_-^{\top}RU_-+X_{+}^{\top}P_GX_{+})Z + P_G-Q \\
	&\preceq \left(\left(\|U_-\|^2\|R\|+\|X_{+}\|^2J(G)+J(G)\right)\frac{1}{\underline{\sigma}(Q)}-1\right)Q\\
	&:= \xi Q.
	\end{align*}
	

	Then, it follows from the definition of $P_G'[Z]$ that
	$$
	P_G'[Z] = \sum_{i=0}^{\infty} (G^{\top}X_{+}^{\top})^i(Z^{\top}E_G + E_G^{\top}Z)(X_{+}G)^i \preceq \xi P_G,
	$$
	and hence $\|P_G'[Z]\|_F \leq \xi J(G)$. Finally, we can bound the Hessian by
	\begin{align*}
	\|\nabla^2 J(G)\| \leq 2\|U_-\|^2\|R\| \frac{J(G)}{\underline{\sigma}(Q)} + (\xi+2) \|X_{+}\|_F^2 \frac{J^2(G)}{\underline{\sigma}(Q)}.
	\end{align*}
	
	Noting that $J(G)\leq a$, the proof is completed.
\end{appendices}

\end{document}